  \documentclass[twoside,reqno,11pt]{fcaa}  

\usepackage{graphicx}
\usepackage{epsfig}
\usepackage{amsthm}
\usepackage{amsmath}
\usepackage{latexsym}
\usepackage{amsfonts}
\usepackage{amssymb}

\newtheoremstyle{theorem}
  {15pt}          
  {15pt}  
  {\sl}  
  {\parindent}
  {\sc}  
  {. }   
  { }    
  {}     
\theoremstyle{theorem}

\newtheoremstyle{defi}
  {15pt}          
  {15pt}  
  {\rm}  
  {\parindent}     
  {\sc}  
  {. }    
  { }    
  {}     
\theoremstyle{defi}

 \def\proofname{\indent\rm P\ r\ o\ o\ f.}
 \def\proofend{\hfill$\Box$}
 \def\qed{\hfill$\Box$} 

 \usepackage{hyperref} 

 \def\theequation{\arabic{section}.\arabic{equation}}

 \def\themycopyrightfootnote{\vspace*{3pt}
 \copyright \, Year\,  Diogenes  Co., Sofia
   \par  \noindent pp. xxx--xxx, DOI: ......................
   \hfill  \vspace*{-36pt}
  }

 \title[ON REPRESENTATION FORMULAS FOR SOLUTIONS \dots]
 {ON REPRESENTATION FORMULAS FOR SOLUTIONS OF LINEAR DIFFERENTIAL EQUATIONS WITH CAPUTO FRACTIONAL DERIVATIVES \\ [3pt] IN ``FCAA'' JOURNAL}
 \author[\normalsize M. I. Gomoyunov]{\normalsize Mikhail I. Gomoyunov}

 \newtheorem{theo}{Theorem}[section]
 \newtheorem{lem}{Lemma}[section]
 \newtheorem{prop}{Proposition}[section]
 \newtheorem{cor}{Corollary}[section]

 \def\Id{\operatorname{Id}_n}
 \def\Linf{\operatorname{L}^\infty}
 \def\ess{\operatorname{ess}}
 \def\C{\operatorname{C}}
 \def\AC{\operatorname{AC}^{\, \alpha}}
 \def\I{{I}}
 \def\D{D}
 \def\R{{R}}
 \def\K{{K}}
 \def\J{{J}}
 \def\k{{k}}
 \def\T{{\mathcal{A}}}
 \def\rd{\mathrm{d}}
 
\begin{document}

 \vbox to 2.5cm { \vfill }


 \bigskip \medskip

 \begin{abstract}

    In the paper, a linear differential equation with variable coefficients and a Caputo fractional derivative is considered.
    For this equation, a Cauchy problem is studied, when an initial condition is given at an intermediate point that does not necessarily coincide with the initial point of the fractional differential operator.
    A detailed analysis of basic properties of the fundamental solution matrix is carried out.
    In particular, the H\"{o}lder continuity of this matrix with respect to both variables is proved, and its dual definition is given.
    Based on this, two representation formulas for the solution of the Cauchy problem are proposed and justified.

 \medskip

{\it MSC 2010\/}: Primary 26A33;
                  Secondary 34A08, 34A30

 \smallskip

{\it Key Words and Phrases}:
    linear fractional differential equation,
    fundamental solution matrix,
    representation formula
    
 \end{abstract}

 \maketitle

 \vspace*{-16pt}


\section*{Introduction}\label{Sec:1}

Representation formulas for solutions of linear differential equations play an important role in differential equations theory.
In particular, they constitute the basis for the study and construction of solutions of various problems arising in stability, control, and differential games theories for linear dynamical systems.

For many classes of linear differential equations with fractional order derivatives, such representation formulas are known.
For various linear equations with constant coefficients, they are given {\it e.g.} in
\cite{Atanackovic_Dolicanin_Pilipovic_Stankovic_2014},
\cite{Bonilla_Rivero_Trujillo_2007},
\cite{Chikriy_Matichin_2008},
\cite[Section~7.1]{Diethelm_2010},
\cite{Duan_2018},
\cite{Idczak_Kamocki_2011},
\cite[Chapter~5]{Kilbas_Srivastava_Trujillo_2006},
\cite{Pskhu_2011}.
The case of variable coefficients is significantly less investigated, see {\it e.g.}
\cite{Bourdin_2018},
\cite[Sect.~7.1]{Diethelm_2010},
\cite{Kaczorek_Idczak_2017}.
Let us particularly note paper \cite{Bourdin_2018}, where not only representation formulas for solutions are derived, but also fundamental solution matrices are introduced and studied in detail under sufficiently general assumptions and for both Riemann--Liouville and Caputo fractional derivatives.
Let us also note that, {\it e.g.} in \cite{Zahariev_Kiskinov_2018}, \cite{Zhang_Wu_2014}, representation formulas are proposed for solutions of fractional order linear functional-differential equations with delays.

However, in all these works, in the case of Caputo derivatives, differential equations are considered only with initial conditions given by a value $x(t_0)$ of a solution at a point $t_0$ that coincides with the initial point of the fractional differential operator.
The goal of this paper is to propose representation formulas for solutions in the case of general initial conditions.
Due to a non-local structure of fractional differential operators, as an initial data at an intermediate point $t_\ast,$ it is natural to consider not only a value $x(t_\ast)$ of a solution at this point, but all the values $x(t)$ for $t \in [t_0, t_\ast],$ see {\it e.g.} \cite{Cong_Tuan_2017}, \cite{Gomoyunov_2019_DGAA}.

The necessity of representation formulas for solutions of linear differential equations with such intermediate initial conditions arises, for example, in the study of control problems and differential games.
In particular, these formulas play a key role in constructing effective and numerically realizable optimal feedback control strategies.
For one of the approaches, successfully developed for various classes of functional-differential equations with delays, the reader is referred to
\cite{Gomoyunov_Lukoyanov_2012},
\cite{Gomoyunov_Lukoyanov_2018_PSIM},
\cite{Krasovskii_Krasovskii_1995},
\cite{Krasovskii_Lukoyanov_1996},
\cite{Lukoyanov_Reshetova_1998}
(see also \cite{Lukoyanov_Gomoyunov_2018_DGAA}).
Bearing in mind such further applications, let us note that the choice of admissible initial functions at intermediate initial points should naturally agree with a notion of a solution of a differential equation.

The rest of the paper is organized as follows.
In Section~\ref{Section_Preliminary}, the definitions and basic properties of fractional order integrals and derivatives are recalled, and also some preliminary results are proved.
In Section~\ref{Section_Cauchy_problem}, a linear differential equation with variable coefficients and a Caputo fractional derivative is considered.
For this equation, a Cauchy problem is formulated, when an initial condition is given at an intermediate point.
A notion of a solution of this Cauchy problem is introduced, a statement on the existence and uniqueness of such a solution is given.
In Section~\ref{Section_Fundamental_matrix}, the fundamental solution matrix is considered, and its properties are studied, including the H\"{o}lder continuity with respect to both variables and a dual definition.
A comparison with the corresponding results of \cite{Bourdin_2018} is carried out.
In Section~\ref{Section_Representation_formulas}, two representation formulas for the solution of the considered Cauchy problem are proposed and justified.

\section{Notations, definitions, and preliminary results}
\label{Section_Preliminary}

\setcounter{section}{1}
\setcounter{equation}{0}\setcounter{theorem}{0}

The basic notations used in the paper are standard.
Let a number $n \in \mathbb{N}$ be fixed.
Let $\mathbb{R}^n$ and $\mathbb{R}^{n \times n}$ be the spaces of $n$-dimensional vectors and $(n \times n)$-matrices, and let $\Id \in \mathbb{R}^{n \times n}$ stand for the identity matrix.
By $\|\cdot\|,$ we denote a norm in $\mathbb{R}^n$ and the corresponding norm in $\mathbb{R}^{n \times n}.$

Let numbers $a, b \in \mathbb{R},$ $a \leq b,$ be fixed, and let $X$ be one of the spaces $\mathbb{R}^n$ or $\mathbb{R}^{n \times n}.$
Let us consider a function $\varphi: [a, b] \rightarrow X,$ for which we also use the notation $\varphi(\cdot).$
We say that $\varphi(\cdot) \in \Linf([a, b], X),$ if the function $\varphi(\cdot)$ is (Lebesgue) measurable and $\|\varphi(\cdot)\|_{[a, b]} < \infty,$ where we denote $\|\varphi(\cdot)\|_{[a, b]} = \ess \sup_{t \in [a, b]} \|\varphi(t)\|$ for $a < b$ and $\|\varphi(\cdot)\|_{[a, b]} = \|\varphi(a)\|$ for $a = b.$
Further, we say that $\varphi(\cdot) \in \C([a, b], X),$ if the function $\varphi(\cdot)$ is continuous.
Let us note that, in the case when $a = b,$ the spaces $\Linf([a, b], X)$ and $\C([a, b], X)$ can be identified with the space $X.$

\subsection{Fractional integrals and derivatives}

Let a number $\alpha \in (0, 1)$ be fixed.
For a function $\varphi: [a, b] \rightarrow X,$ the left-sided and right-sided Riemann--Liouville fractional integrals of the order $\alpha$ are defined for $t \in [a, b]$ by
\begin{equation*}
    (\I^\alpha_{a +} \varphi) (t)
    = \frac{1}{\Gamma(\alpha)} \int_{a}^{t} \frac{\varphi(\tau)}{(t - \tau)^{1 - \alpha}} \, \rd \tau,
    \ \ (\I^\alpha_{b -} \varphi) (t)
    = \frac{1}{\Gamma(\alpha)} \int_{t}^{b} \frac{\varphi(\tau)}{(\tau - t)^{1 - \alpha}} \, \rd \tau,
\end{equation*}
where $\Gamma$ is the gamma-function.
According to \cite[Proposition~2.1]{Gomoyunov_2018_FCAA} (see also the references therein), the following proposition holds.
\begin{prop} \label{Prop_I}
    For any $\varphi(\cdot) \in \Linf([a, b], X),$ the values $(\I^\alpha_{a +} \varphi) (t)$ and $(\I^\alpha_{b -} \varphi) (t)$ are correctly defined for every $t \in [a, b].$
    The inequalities
    \begin{equation*}
        \begin{array}{c}
            \|(\I_{a +}^\alpha \varphi)(t_1) - (\I_{a +}^\alpha \varphi)(t_2)\|
            \leq H_\I \|\varphi(\cdot)\|_{[a, b]} |t_1 - t_2|^\alpha, \\[0.5em]
            \|(\I_{b -}^\alpha \varphi)(t_1) - (\I_{b -}^\alpha \varphi)(t_2)\|
            \leq H_\I \|\varphi(\cdot)\|_{[a, b]} |t_1 - t_2|^\alpha
        \end{array}
    \end{equation*}
    are valid for $t_1, t_2 \in [a, b]$ with the number $H_\I = 2 / \Gamma(\alpha + 1).$
    In particular, the inclusions $(\I_{a +}^\alpha \varphi)(\cdot), (\I_{b -}^\alpha \varphi)(\cdot) \in \C([a, b], X)$ hold.
\end{prop}

For a function $x: [a, b] \rightarrow X,$ the left-sided Riemann--Liouville and Caputo fractional derivatives of the order $\alpha$ are defined for $t \in [a, b]$ by
\begin{equation*}
    \begin{array}{c}
        \displaystyle
        (\D_{a +}^{\, \alpha} x)(t)
        = \frac{\rd}{\rd t} (\I_{a +}^{1 - \alpha} x)(t)
        = \frac{1}{\Gamma(1 - \alpha)} \frac{\rd}{\rd t} \int_{a}^{t} \frac{x(\tau)}{(t - \tau)^\alpha} \, \rd \tau, \\[1em]
        \displaystyle
        (^C \D_{a +}^{\, \alpha} x)(t)
        = \frac{\rd}{\rd t} \big(\I_{a +}^{1 - \alpha} (x(\cdot) - x(a)) \big) (t)
        = \frac{1}{\Gamma(1 - \alpha)} \frac{\rd}{\rd t} \int_{a}^{t} \frac{x(\tau) - x(a)}{(t - \tau)^\alpha} \, \rd \tau.
    \end{array}
\end{equation*}
We say that $x(\cdot) \in \AC([a, b], X),$ if there exists $\varphi(\cdot) \in \Linf([a, b], X)$ such that $x(t) = x(a) + (\I_{a +}^\alpha \varphi) (t),$ $t \in [a, b].$
Let us note that, if $a = b,$ then $\AC([a, b], X)$ can be identified with $X.$
According to \cite[Proposition~2.2]{Gomoyunov_2018_FCAA} (see also the references therein), the following proposition is valid.
\begin{prop} \label{Proposition_D}
    For any $x(\cdot) \in \AC([a, b], X),$ the value $(^C \D_{a +}^{\, \alpha} x)(t)$ is correctly defined for almost every $t \in [a, b].$
    Moreover, the inclusion $(^C \D_{a +}^{\, \alpha} x)(\cdot) \in \Linf([a, b], X)$ holds (i.e., there exists $\varphi(\cdot) \in \Linf([a, b], X)$ such that $\varphi(t) = (^C \D_{a +}^{\, \alpha} x)(t)$ for almost every $t \in [a, b]$), and
    \begin{equation*}
        \big(\I_{a +}^\alpha (^C \D_{a +}^{\, \alpha} x) \big)(t)
        = x(t) - x(a),
        \quad t \in [a, b].
    \end{equation*}
\end{prop}

\subsection{Auxiliary operators $\J^\alpha_{a +}$ and $\J^\alpha_{b -}$}

For a function $\varphi: [a, b] \rightarrow X,$ let us denote
\begin{equation*}
    \begin{array}{c}
        \displaystyle
        (\R_{a +}^{\, \alpha} \varphi)(t)
        = \frac{(1 - \alpha) \sin(\alpha \pi)}{\pi} \int_{a}^{t} \K(t - a, \tau - a) \varphi(\tau) \, \rd \tau, \\[1em]
        \displaystyle
        (\R_{b -}^{\, \alpha} \varphi)(t)
        = \frac{(1 - \alpha) \sin(\alpha \pi)}{\pi} \int_{t}^{b} \K(b - t, b - \tau) \varphi(\tau) \, \rd \tau,
    \end{array}
\end{equation*}
where $t \in [a, b]$ and, for $(\xi, \tau) \in (0, \infty) \times (0, \infty)$ such that $\xi > \tau,$
\begin{equation*}
    \K(\xi, \tau)
    = \frac{1}{\tau^{1 - \alpha}} \int_{0}^{1} \frac{\eta^\alpha}{(1 - \eta)^\alpha (\tau + \eta (\xi - \tau))^\alpha} \, \rd \eta.
\end{equation*}
One can show that the function $\K$ is correctly defined, non-negative, continuous, and homogenous of degree $-1.$
Furthermore,
\begin{equation*}
    \K(1, \tau)
    \leq \frac{1}{(1 - \alpha) \tau^{1 - \alpha}},
    \quad \tau \in (0, 1).
\end{equation*}
Based on these properties of $\K,$ one can obtain the result below.
\begin{prop} \label{Proposition_R}
    For any $\varphi(\cdot) \in \Linf([a, b], X),$ the values $(\R_{a +}^{\, \alpha} \varphi)(t)$ and $(\R_{b -}^{\, \alpha} \varphi)(t)$ are correctly defined for every $t \in [a, b].$
    Moreover, the functions $(\R_{a +}^{\, \alpha} \varphi)(\cdot)$ and $(\R_{b -}^{\, \alpha} \varphi)(\cdot)$ are measurable, and the inequalities
    \begin{equation*}
        \|(\R_{a +}^{\, \alpha} \varphi)(t)\|
        \leq M_\R \|\varphi(\cdot)\|_{[a, t]},
        \quad
        \|(\R_{b -}^{\, \alpha} \varphi)(t)\|
        \leq M_\R \|\varphi(\cdot)\|_{[t, b]},
        \quad t \in [a, b],
    \end{equation*}
    are valid with the number $M_\R = (\alpha \pi)^{-1} \sin(\alpha \pi).$
    In particular, the inclusions $(\R_{a +}^{\, \alpha} \varphi)(\cdot), (\R_{b -}^{\, \alpha} \varphi)(\cdot) \in \Linf([a, b], X)$ hold.
\end{prop}

For a function $\varphi: [a, b] \rightarrow X,$ let us denote
\begin{equation*}
    \begin{array}{c}
        \displaystyle
        (\J_{a +}^\alpha \varphi)(t)
        = \frac{(t - a)^{1 - \alpha}}{\Gamma(\alpha)}
        \int_{a}^{t} \frac{\varphi(\tau)}{(t - \tau)^{1 - \alpha} (\tau - a)^{1 - \alpha}} \, \rd \tau, \\[1em]
        \displaystyle
        (\J_{b -}^\alpha \varphi)(t)
        = \frac{(b - t)^{1 - \alpha}}{\Gamma(\alpha)}
        \int_{t}^{b} \frac{\varphi(\tau)}{(b - \tau)^{1 - \alpha} (\tau - t)^{1 - \alpha}} \, \rd \tau,
    \end{array}
\end{equation*}
where $t \in [a, b]$.
By the scheme from \cite[Lemma~3.1]{Samko_Kilbas_Marichev_1993}, one can prove the following lemma.
\begin{lem} \label{Lemma_J}
    For any $\varphi(\cdot) \in \Linf([a, b], X),$ the values $(\J_{a +}^\alpha \varphi)(t)$ and $(\J_{b -}^\alpha \varphi)(t)$ are correctly defined for every $t \in [a, b],$ and the equalities
    \begin{equation} \label{J_R}
        \begin{array}{c}
            \displaystyle
            (\J_{a +}^\alpha \varphi)(t)
            = \frac{1}{\Gamma(\alpha)} \int_{a}^{t} \frac{\varphi(\tau) + (\R_{a +}^{\, \alpha} \varphi)(\tau)}{(t - \tau)^{1 - \alpha}} \, \rd \tau, \\[1em]
            \displaystyle
            (\J_{b -}^\alpha \varphi)(t)
            = \frac{1}{\Gamma(\alpha)} \int_{t}^{b} \frac{\varphi(\tau) + (\R_{b -}^{\, \alpha} \varphi)(\tau)}{(\tau - t)^{1 - \alpha}} \, \rd \tau
        \end{array}
    \end{equation}
    are valid for $t \in [a, b].$
\end{lem}

\begin{cor} \label{Corollary_J_Holder}
    For any $\varphi(\cdot) \in \Linf([a, b], X),$ the inequalities
    \begin{equation} \label{J_Holder}
        \begin{array}{c}
            \|(\J_{a +}^\alpha \varphi)(t_1) - (\J_{a +}^\alpha \varphi)(t_2)\|
            \leq H_\J \|\varphi(\cdot)\|_{[a, b]} |t_1 - t_2|^\alpha, \\[0.5em]
            \|(\J_{b -}^\alpha \varphi)(t_1) - (\J_{b -}^\alpha \varphi)(t_2)\|
            \leq H_\J \|\varphi(\cdot)\|_{[a, b]} |t_1 - t_2|^\alpha,
        \end{array}
    \end{equation}
    where $t_1, t_2 \in [a, b],$ and
    \begin{equation} \label{J_Boundedness}
        \begin{array}{c}
            \displaystyle
            \| (\J_{a +}^\alpha \varphi)(t) \|
            \leq \frac{M_\J}{\Gamma(\alpha)} \int_{a}^{t} \frac{\|\varphi(\cdot)\|_{[a, \tau]}}{(t - \tau)^{1 - \alpha}} \, \rd \tau, \\[1em]
            \displaystyle
            \| (\J_{b -}^\alpha \varphi)(t) \|
            \leq \frac{M_\J}{\Gamma(\alpha)} \int_{t}^{b} \frac{\|\varphi(\cdot)\|_{[\tau, b]}}{(\tau - t)^{1 - \alpha}} \, \rd \tau,
        \end{array}
    \end{equation}
    where $t \in [a, b],$ are valid with the numbers
    \begin{equation*}
        H_\J
        = (1 + M_\R) H_\I
        = \frac{2}{\Gamma(\alpha + 1)} \Big( 1 + \frac{\sin(\alpha \pi)}{\alpha \pi} \Big),
        \ \ M_\J
        = 1 + M_\R
        = 1 + \frac{\sin(\alpha \pi)}{\alpha \pi}.
    \end{equation*}
    In particular, the inclusions $(\J_{a +}^\alpha \varphi)(\cdot), (\J_{b -}^\alpha \varphi)(\cdot) \in \C([a, b], X)$ hold.
\end{cor}
\proof
    Let $\varphi(\cdot) \in \Linf([a, b], X)$ be fixed.
    According to the first equality in (\ref{J_R}), for the function $\psi(t) = \varphi(t) + (\R_{a +}^{\, \alpha} \varphi)(t),$ $t \in [a, b],$ we have $(\J_{a +}^\alpha \varphi)(t) = (\I_{a +}^\alpha \psi)(t),$ $t \in [a, b].$
    By Proposition~\ref{Proposition_R},
    \begin{equation*}
        \|\psi(t)\|
        \leq \|\varphi(t)\| + \|(\R_{a +}^{\, \alpha} \varphi)(t)\|
        \leq (1 + M_\R) \|\varphi(\cdot)\|_{[a, t]}
        \text{ for a.e. } t \in [a, b].
    \end{equation*}
    Hence, for $t \in [a, b]$,
    \begin{multline*}
        \|(\J_{a +}^\alpha \varphi)(t)\|
        = \|(\I_{a +}^\alpha \psi)(t)\| \\
        \leq \frac{1}{\Gamma(\alpha)} \int_{a}^{t} \frac{\|\psi(\tau)\|}{(t - \tau)^{1 - \alpha}} \, \rd \tau
        \leq \frac{1 + M_\R}{\Gamma(\alpha)}
        \int_{a}^t \frac{\|\varphi(\cdot)\|_{[a, \tau]}}{(t - \tau)^{1 - \alpha}} \, \rd \tau,
    \end{multline*}
    and, due to Proposition~\ref{Prop_I}, for $t_1, t_2 \in [a, b],$
    \begin{multline*}
        \|(\J_{a +}^\alpha \varphi)(t_1) - (\J_{a +}^\alpha \varphi)(t_2)\|
        = \|(\I_{a +}^\alpha \psi)(t_1) - (\J_{a +}^\alpha \psi)(t_2)\| \\
        \leq H_\I \|\psi(\cdot)\|_{[a, b]} |t_1 - t_2|^\alpha
        \leq (1 + M_\R) H_\I \|\varphi(\cdot)\|_{[a, b]} |t_1 - t_2|^\alpha.
    \end{multline*}
    Thus, the first inequalities in (\ref{J_Holder}) and (\ref{J_Boundedness}) are proved.
    The validity of the second inequalities in (\ref{J_Holder}) and (\ref{J_Boundedness}) can be shown in a similar way.
\proofend

\section{Linear differential equations with Caputo derivatives}
\label{Section_Cauchy_problem}

\setcounter{section}{2}
\setcounter{equation}{0}\setcounter{theorem}{0}

Let us fix numbers $t_0, \vartheta \in \mathbb{R},$ $t_0 < \vartheta,$ and functions
\begin{equation} \label{A_b}
    A(\cdot) \in \Linf([t_0, \vartheta], \mathbb{R}^{n \times n}),
    \quad b(\cdot) \in \Linf([t_0, \vartheta], \mathbb{R}^n).
\end{equation}
For a number $t_\ast \in [t_0, \vartheta)$ and a function $w_\ast(\cdot) \in \AC([t_0, t_\ast], \mathbb{R}^n),$ let us consider the following Cauchy problem for the linear differential equation with the Caputo fractional derivative
\begin{equation} \label{differential_equation}
    ({}^C \D_{t_0 +}^{\, \alpha} x) (t) = A(t) x(t) + b(t),
    \quad x(t) \in \mathbb{R}^n, \quad t \in [t_\ast, \vartheta],
\end{equation}
and the initial condition
\begin{equation} \label{initial_condition}
    x(t) = w_\ast(t), \quad t \in [t_0, t_\ast].
\end{equation}
A function $x: [t_0, \vartheta] \rightarrow \mathbb{R}^n$ is called a solution of Cauchy problem (\ref{differential_equation}), (\ref{initial_condition}), if the inclusion $x(\cdot) \in \AC([t_0, \vartheta], \mathbb{R}^n)$ is valid, initial condition (\ref{initial_condition}) holds, and differential equation (\ref{differential_equation}) is satisfied for almost every $t \in [t_\ast, \vartheta].$
By the scheme from \cite[Proposition~2]{Gomoyunov_2019_DGAA} (see also \cite[Theorem~3.1]{Gomoyunov_2018_FCAA}), one can prove the proposition below.
\begin{prop} \label{Proposition_x}
    For any $t_\ast \in [t_0, \vartheta)$ and any $w_\ast(\cdot) \in \AC([t_0, t_\ast], \mathbb{R}^n),$ there exists a unique solution of Cauchy problem (\ref{differential_equation}), (\ref{initial_condition}).
    This solution is a unique function $x(\cdot) \in \C([t_0, \vartheta], \mathbb{R}^n)$ that satisfies initial condition (\ref{initial_condition}) and the integral equation
    \begin{multline} \label{integral_equation}
        x(t) = w_\ast(t_0)
        + \frac{1}{\Gamma(\alpha)} \int_{t_0}^{t_\ast} \frac{(^C \D_{t_0 +}^{\, \alpha} w_\ast)(\tau)}{(t - \tau)^{1 - \alpha}} \, \rd \tau \\
        + \frac{1}{\Gamma(\alpha)} \int_{t_\ast}^{t} \frac{A(\tau) x(\tau) + b(\tau)}{(t - \tau)^{1 - \alpha}} \, \rd \tau,
        \quad t \in [t_\ast, \vartheta].
    \end{multline}
\end{prop}

The goal of the paper is to obtain representation formulas for the solution $x(\cdot)$ of Cauchy problem (\ref{differential_equation}), (\ref{initial_condition}).
These formulas are based on an appropriate notion of the fundamental solution matrix, which is introduced and studied in the next section.

\section{The fundamental solution matrix}
\label{Section_Fundamental_matrix}

\setcounter{section}{3}
\setcounter{equation}{0}\setcounter{theorem}{0}

Let us denote
\begin{equation*}
    \Omega
    = \big\{ (t, s) \in [t_0, \vartheta] \times [t_0, \vartheta]: \, t \geq s \big\},
    \quad
    \Omega_0
    = \big\{ (t, s) \in \Omega: \, t > s \big\}.
\end{equation*}
Following \cite{Bourdin_2018}, one can consider the function $\Omega_0 \ni (t, s) \mapsto Z(t, s) \in \mathbb{R}^{n \times n}$ that is defined for any $s \in [t_0, \vartheta)$ as a solution of the Cauchy problem for the linear homogenous differential equation with the Riemann--Liuoville derivative
\begin{equation} \label{differential_equation_Z}
    (\D_{s +}^{\, \alpha} Z(\cdot, s))(t)
    = A(t) Z(t, s),
    \quad t \in (s, \vartheta],
\end{equation}
and the initial condition
\begin{equation} \label{initial_condition_Z}
    (\I_{s +}^{1 - \alpha} Z(\cdot, s)) (s)
    = \Id,
\end{equation}
or, which is equivalent, as a solution of the corresponding integral equation
\begin{equation}\label{integral_equation_Z}
    Z(t, s)
    = \frac{\Id}{\Gamma(\alpha) (t - s)^{1 - \alpha}}
    + \frac{1}{\Gamma(\alpha)} \int_{s}^{t} \frac{A(\tau) Z(\tau, s)}{(t - \tau)^{1 - \alpha}} \, \rd \tau,
    \quad t \in (s, \vartheta].
\end{equation}
In \cite{Bourdin_2018}, a detailed analysis of Cauchy problem (\ref{differential_equation_Z}), (\ref{initial_condition_Z}) and integral equation (\ref{integral_equation_Z}) is carried out, and, in terms of the function $Z,$ a representation formula for the solution $x(\cdot)$ of the original Cauchy problem (\ref{differential_equation}), (\ref{initial_condition}) is given in a particular case when $t_\ast = t_0.$
Let us note that the function $Z$ may have a singularity when $t$ goes to $s,$ which significantly complicates the study and use of this function.
For example, in the case of constant coefficients, i.e., when $A(t) = A_0 \in \mathbb{R}^{n \times n},$ $t \in [t_0, \vartheta],$ we have (see \cite[Example~2]{Bourdin_2018})
\begin{equation*}
    Z(t, s)
    = \frac{E_{\alpha, \alpha} ((t - s)^{\alpha} A_0)}{(t - s)^{1 - \alpha}},
    \quad (t, s) \in \Omega_0,
\end{equation*}
where $E_{\alpha, \alpha}$ is the two-parametric Mittag-Leffler function.

One of the questions here is related to derivation of the properties of the function $[t_0, t) \ni s \mapsto Z(t, s) \in \mathbb{R}^{n \times n}$ for a fixed $t \in (t_0, \vartheta]$ directly from relation (\ref{integral_equation_Z}).
In \cite[Lemma~5]{Bourdin_2018}, the measurability of this function is stated, but this fact is not discussed in the given proof.
Nevertheless, it plays a key role in the proofs of the representation formula and duality theorem \cite[Theorems~5 and~7]{Bourdin_2018}.
Since this question is also important for the present paper, it is suggested to pay more attention to it.
In this way, it is proposed to introduce the function $(t - s)^{1 - \alpha} Z(t, s),$ $(t, s) \in \Omega,$ and study the properties of this function instead of $Z$ itself.

Namely, let us consider the function $\Omega \ni (t, s) \mapsto F(t, s) \in \mathbb{R}^{n \times n}$ that is defined for any $s \in [t_0, \vartheta]$ as a continuous solution of the integral equation
\begin{equation}\label{F}
    F(t, s)
    = \frac{\Id}{\Gamma(\alpha)}
    + \frac{(t - s)^{1 - \alpha}}{\Gamma(\alpha)} \int_{s}^{t} \frac{A(\tau) F(\tau, s)}{(t - \tau)^{1 - \alpha} (\tau - s)^{1 - \alpha}} \, \rd \tau,
    \quad t \in [s, \vartheta].
\end{equation}
The propositions below establish basic properties of the function $F.$
Despite the fact that some of them can be derived from the corresponding properties of the function $Z$ and the equality
\begin{equation} \label{F_Z}
    F(t, s)
    = (t - s)^{1 - \alpha} Z(t, s),
    \quad (t, s) \in \Omega_0,
\end{equation}
their independent and complete proofs are given, since the constructions and technique used are slightly different from \cite{Bourdin_2018}.

\begin{prop} \label{Proposition_Existence_and_Boundedness_F}
    For any $s \in [t_0, \vartheta],$ there exists a unique continuous function $[s, \vartheta] \ni t \mapsto F(t, s) \in \mathbb{R}^{n \times n}$ satisfying integral equation (\ref{F}).
    Moreover, there exists a number $M_F > 0$ such that
    \begin{equation} \label{M_F}
        \|F(t, s)\| \leq M_F,
        \quad (t, s) \in \Omega.
    \end{equation}
\end{prop}
\proof
    For $s = \vartheta,$ we obtain from (\ref{F}) that $F(\vartheta, \vartheta) = \Id / \Gamma(\alpha).$
    Let $s \in [t_0, \vartheta)$ be fixed.
    Let us denote
    \begin{equation} \label{M_A}
        M_A
        = \|A(\cdot)\|_{[t_0, \vartheta]}
    \end{equation}
    and choose a number $\k > 0$ such that $\k^{- \alpha} M_A M_\J < 1,$ where $M_\J$ is the number from Corollary~\ref{Corollary_J_Holder}.
    Let us consider the Banach space $\C_e([s, \vartheta], \mathbb{R}^{n \times n})$ of functions $\varphi(\cdot) \in \C([s, \vartheta], \mathbb{R}^{n \times n})$ with the Bielecki norm
    \begin{equation*}
        \|\varphi(\cdot)\|_e
        = \max_{t \in [s, \vartheta]} (\|\varphi(t)\|e^{- (t - s) \k}).
    \end{equation*}
    By the right-hand side of integral equation (\ref{F}), let us define the operator $\T: \C_e([s, \vartheta], \mathbb{R}^{n \times n}) \rightarrow \C_e([s, \vartheta], \mathbb{R}^{n \times n})$ as follows:
    \begin{equation*}
        \begin{array}{c}
            (\T \varphi)(t)
            = \displaystyle \frac{\Id}{\Gamma(\alpha)}
            + \frac{(t - s)^{1 - \alpha}}{\Gamma(\alpha)} \int_{s}^{t} \frac{A(\tau) \varphi(\tau)}{(t - \tau)^{1 - \alpha} (\tau - s)^{1 - \alpha}} \, \rd \tau, \\[1em]
            t \in [s, \vartheta],
            \quad \varphi(\cdot) \in \C_e([s, \vartheta], \mathbb{R}^{n \times n}).
        \end{array}
    \end{equation*}
    Let us note that, for any $\varphi(\cdot) \in \C_e([s, \vartheta], \mathbb{R}^{n \times n}),$
    \begin{equation*}
        (\T \varphi)(t)
        = \frac{\Id}{\Gamma(\alpha)}
        + (\J_{s +}^\alpha \psi)(t),
        \quad t \in [s, \vartheta],
    \end{equation*}
    where $\psi(t) = A(t) \varphi(t),$ $t \in [s, \vartheta].$
    According to the first inclusion in (\ref{A_b}), we have $\psi(\cdot) \in \Linf([s, \vartheta], \mathbb{R}^{n \times n}),$ and, therefore, it follows from Corollary~\ref{Corollary_J_Holder} that $(\T \varphi)(\cdot) \in \C_e([s, \vartheta], \mathbb{R}^{n \times n}).$
    Hence, the operator $\T$ is correctly defined.
    Let us show that this operator is a contraction.
    Let us fix $\varphi_1(\cdot), \varphi_2(\cdot) \in \C_e([s, \vartheta], \mathbb{R}^{n \times n}).$
    For any $t \in [s, \vartheta],$ we derive
    \begin{equation*}
        \|(\T \varphi_1)(t) - (\T \varphi_2)(t)\|
        \leq \frac{(t - s)^{1 - \alpha} M_A}{\Gamma(\alpha)}
        \int_{s}^{t} \frac{\|\varphi_1(\tau) - \varphi_2(\tau)\|}{(t - \tau)^{1 - \alpha} (\tau - s)^{1 - \alpha}} \, \rd \tau.
    \end{equation*}
    Since
    \begin{multline*}
        \max_{\xi \in [s, \tau]} \|\varphi_1(\xi) - \varphi_2(\xi)\|
        \leq \max_{\xi \in [s, \tau]} (\| \varphi_1(\xi) - \varphi_2(\xi) \|
        e^{- (\xi - s) \k}) e^{(\tau - s) \k} \\
        \leq \| \varphi_1(\cdot) - \varphi_2(\cdot) \|_e e^{(\tau - s) \k},
        \quad \tau \in [s, \vartheta],
    \end{multline*}
    then, due to Corollary~\ref{Corollary_J_Holder}, we obtain
    \begin{multline*}
        \frac{(t - s)^{1 - \alpha}}{\Gamma(\alpha)}
        \int_{s}^{t} \frac{\|\varphi_1(\tau) - \varphi_2(\tau)\|}{(t - \tau)^{1 - \alpha} (\tau - s)^{1 - \alpha}} \, \rd \tau \\
        \leq \frac{M_\J \| \varphi_1(\cdot) - \varphi_2(\cdot) \|_e}{\Gamma(\alpha)}
        \int_{s}^{t} \frac{ e^{(\tau - s) \k}}{(t - \tau)^{1 - \alpha}} \, \rd \tau.
    \end{multline*}
    Consequently, taking the inequality
    \begin{equation*}
        \frac{1}{\Gamma(\alpha)} \int_{s}^{t} \frac{e^{(\tau - s) \k}}{(t - \tau)^{1 - \alpha}} \, \rd \tau
        \leq e^{(t - s) \k} \k^{- \alpha}
    \end{equation*}
    into account (see {\it e.g.} the proof of \cite[Theorem~1]{Bourdin_2018}), we conclude
    \begin{equation} \label{Proposition_Existence_and_Boundedness_F_p_contraction}
        \|(\T \varphi_1)(\cdot) - (\T \varphi_2)(\cdot)\|_e
        \leq \k^{- \alpha} M_A M_\J \| \varphi_1(\cdot) - \varphi_2(\cdot) \|_e.
    \end{equation}
    Thus, due to the choice of the number $\k,$ the operator $\T$ is a contraction.
    By the Banach contraction principle, this operator has a unique fixed point, which is a unique continuous solution $[s, \vartheta] \ni t \mapsto F(t, s) \in \mathbb{R}^{n \times n}$ of (\ref{F}).

    Further, let us prove that estimate (\ref{M_F}) is valid with the number
    \begin{equation*}
        M_F
        = \frac{e^{(\vartheta - t_0) \k}}{\Gamma(\alpha) (1 - \k^{- \alpha} M_A M_\J)}.
    \end{equation*}
    For $s = \vartheta,$ we have $\|F(\vartheta, \vartheta)\| = 1 / \Gamma(\alpha) \leq M_F.$
    Let $s \in [t_0, \vartheta)$ be fixed.
    Let us introduce the function $\varphi_0(t) = 0 \in \mathbb{R}^{n \times n},$ $t \in [s, \vartheta].$
    Then, according to inequality (\ref{Proposition_Existence_and_Boundedness_F_p_contraction}), we derive
    \begin{multline*}
        \|F(\cdot, s)\|_e
        = \|(\T F(\cdot, s))(\cdot)\|_e \\
        \leq \|(\T F(\cdot, s))(\cdot) - (\T \varphi_0)(\cdot)\|_e + \|(\T \varphi_0)(\cdot)\|_e \\
        \leq \k^{- \alpha} M_A M_\J \|F(\cdot, s)\|_e + \|(\T \varphi_0)(\cdot)\|_e.
    \end{multline*}
    Hence,
    \begin{equation*}
        \|F(\cdot, s)\|_e
        \leq \frac{\|(\T \varphi_0)(\cdot)\|_e}{1 - \k^{- \alpha} M_A M_\J},
    \end{equation*}
    wherefrom, taking into account that $\|(\T \varphi_0)(\cdot)\|_e = 1 / \Gamma(\alpha),$ we conclude
    \begin{equation*}
        \|F(t, s)\|
        \leq \|F(\cdot, s)\|_e e^{(\vartheta - t_0) k}
        \leq M_F,
        \quad t \in [s, \vartheta].
    \end{equation*}
    The proposition is proved.
\proofend

In the proof of the next proposition, the following version of the Bellman--Gronwall lemma is used, which is a direct corollary of \cite[Lemma~6.19]{Diethelm_2010} and \cite[Proposition~2]{Gomoyunov_2019_PFDA}.
\begin{lem} \label{Lem_BG}
    Let numbers $a, b \in \mathbb{R},$ $a \leq b,$ $\varepsilon_1 \geq 0,$ and $\varepsilon_2 \geq 0$ be fixed.
    Then, for any function $\psi(\cdot) \in \C([a, b], \mathbb{R})$ such that
    \begin{equation*}
        0 \leq
        \psi(t)
        \leq \varepsilon_1 + \frac{\varepsilon_2}{\Gamma(\alpha)}
        \int_{a}^{t} \frac{ \|\psi(\cdot)\|_{[a, \tau]}}{(t - \tau)^{1 - \alpha}} \, \rd \tau,
        \quad t \in [a, b],
    \end{equation*}
    the inequality
    \begin{equation*}
        \psi(t)
        \leq \varepsilon_1 E_\alpha((t - a)^\alpha \varepsilon_2),
        \quad t \in [a, b],
    \end{equation*}
    is valid, where $E_\alpha$ is the one-parametric Mittag-Leffler function.
\end{lem}

\begin{prop} \label{Proposition_F_Holder}
    There exists a number $H_F \geq 0$ such that
    \begin{equation} \label{H_F}
        \begin{array}{c}
            \|F(t_1, s_1) - F(t_2, s_2)\|
            \leq H_F (|t_1 - t_2|^\alpha + |s_1 - s_2|^\alpha), \\[0.5em]
            (t_1, s_1), (t_2, s_2) \in \Omega.
        \end{array}
    \end{equation}
    In particular, the function $\Omega \ni (t, s) \mapsto F(t, s) \in \mathbb{R}^{n \times n}$ is continuous.
\end{prop}
\proof
    Let us prove the H\"{o}lder continuity of the function $F$ with respect to the first variable.
    Let $s \in [t_0, \vartheta]$ be fixed.
    From definition (\ref{F}) of this function, it follows that
    \begin{equation*}
        F(t, s)
        = \frac{\Id}{\Gamma(\alpha)} + (\J^\alpha_{s +} \varphi)(t),
        \quad t \in [s, \vartheta],
    \end{equation*}
    where $\varphi(t) = A(t) F(t, s),$ $t \in [s, \vartheta].$
    Due to the first inclusion in (\ref{A_b}) and Proposition~\ref{Proposition_Existence_and_Boundedness_F}, we have $\varphi(\cdot) \in \Linf([s, \vartheta], \mathbb{R}^{n \times n})$ and $\|\varphi(\cdot)\|_{[s, \vartheta]} \leq M_A M_F,$ where $M_A$ is the number from (\ref{M_A}).
    Hence, by Corollary~\ref{Corollary_J_Holder}, we obtain
    \begin{multline} \label{Proposition_F_Holder_p_t_1_t_2}
        \|F(t_1, s) - F(t_2, s)\|
        = \|(\J^\alpha_{s +} \varphi)(t_1) - (\J^\alpha_{s +} \varphi)(t_2)\| \\
        \leq H_\J \|\varphi(\cdot)\|_{[s, \vartheta]} |t_1 - t_2|^\alpha
        \leq H_\J M_A M_F |t_1 - t_2|^\alpha,
        \quad t_1, t_2 \in [s, \vartheta].
    \end{multline}

    Further, let us prove the H\"{o}lder continuity of $F$ with respect to the second variable.
    Let $s_1, s_2 \in [t_0, \vartheta],$ $s_1 < s_2,$ be fixed.
    According to (\ref{F}),
    \begin{multline} \label{Proposition_F_Holder_p_base}
        F(t, s_1) - F(t, s_2)
        = \frac{(t - s_1)^{1 - \alpha}}{\Gamma(\alpha)}
        \int_{s_1}^{t} \frac{A(\tau) F(\tau, s_1)}{(t - \tau)^{1 - \alpha} (\tau - s_1)^{1 - \alpha}} \, \rd \tau \\
        - \frac{(t - s_2)^{1 - \alpha}}{\Gamma(\alpha)}
        \int_{s_2}^{t} \frac{A(\tau) F(\tau, s_1)}{(t - \tau)^{1 - \alpha} (\tau - s_2)^{1 - \alpha}} \, \rd \tau \\
        + \frac{(t - s_2)^{1 - \alpha}}{\Gamma(\alpha)}
        \int_{s_2}^{t} \frac{A(\tau) (F(\tau, s_1) - F(\tau, s_2))}{(t - \tau)^{1 - \alpha} (\tau - s_2)^{1 - \alpha}} \, \rd \tau,
        \quad t \in [s_2, \vartheta].
    \end{multline}
    Let us denote $\psi(t) = \|F(t, s_1) - F(t, s_2)\|,$ $t \in [s_2, \vartheta],$ and $\chi(t) = A(t) F(t, s_1),$ $t \in [s_1, \vartheta].$
    Then, we have $\psi(\cdot) \in \C([s_2, \vartheta], \mathbb{R}),$ $\chi(\cdot) \in \Linf([s_1, \vartheta], \mathbb{R}^{n \times n}),$ and $\|\chi(\cdot)\|_{[s_1, \vartheta]} \leq M_A M_F.$
    Thus, from equality (\ref{Proposition_F_Holder_p_base}), taking Corollary~\ref{Corollary_J_Holder} into account, we derive
    \begin{multline*}
        \psi(t)
        \leq \|(\J^\alpha_{t -} \chi)(s_1) - (\J^\alpha_{t -} \chi)(s_2)\|
        + M_A (\J^\alpha_{s_2 +} \psi)(t) \\
        \leq H_\J M_A M_F (s_2 - s_1)^\alpha
        + \frac{M_A M_\J}{\Gamma(\alpha)}
        \int_{s_2}^{t} \frac{\|\psi(\cdot)\|_{[s_2, \tau]}}{(t - \tau)^{1 - \alpha}} \, \rd \tau,
        \quad t \in [s_2, \vartheta],
    \end{multline*}
    wherefrom, applying Lemma~\ref{Lem_BG}, we conclude
    \begin{multline} \label{Proposition_F_Holder_p_s_1_s_2}
            \|F(t, s_1) - F(t, s_2)\|
            = \psi(t)
            \leq H_\J M_A M_F (s_2 - s_1)^\alpha
            E_\alpha((t - s_2)^\alpha M_A M_\J) \\
            \leq H_\J M_A M_F
            E_\alpha((\vartheta - t_0)^\alpha M_A M_\J) (s_2 - s_1)^\alpha,
            \quad t \in [s_2, \vartheta].
    \end{multline}

    It follows from estimates (\ref{Proposition_F_Holder_p_t_1_t_2}) and (\ref{Proposition_F_Holder_p_s_1_s_2}) that inequality (\ref{H_F}) is valid with the number
    \begin{equation*}
        H_F
        = H_\J M_A M_F
        E_\alpha((\vartheta - t_0)^\alpha M_A M_\J).
    \end{equation*}
    The proposition is proved.
\proofend

Let us note that the H\"{o}lder continuity of the function $[s, \vartheta] \ni t \mapsto F(t, s) \in \mathbb{R}^{n \times n}$ for a fixed $s \in [t_0, \vartheta]$ can be derived from relation (\ref{F_Z}), connecting the functions $F$ and $Z,$ and the continuity properties of the function $(s, \vartheta] \ni t \mapsto Z(t, s) \in \mathbb{R}^{n \times n}$ as a solution of Cauchy problem (\ref{differential_equation_Z}), (\ref{initial_condition_Z}) (see {\it e.g.} \cite[Theorem~6]{Bourdin_2018_2}).

Following \cite[Theorem~7]{Bourdin_2018}, for any $t \in [t_0, \vartheta],$ let us consider the dual integral equation
\begin{equation} \label{G}
    G(t, s)
    = \frac{\Id}{\Gamma(\alpha)}
    + \frac{(t - s)^{1 - \alpha}}{\Gamma(\alpha)} \int_{s}^{t} \frac{G(t, \tau) A(\tau)}{(t - \tau)^{1 - \alpha} (\tau - s)^{1 - \alpha}} \, \rd \tau,
    \ \ s \in [t_0, t].
\end{equation}

\begin{prop} \label{Proposition_F_G}
    For any $t \in [t_0, \vartheta],$ there exists a unique continuous function $[t_0, t] \ni s \mapsto G(t, s) \in \mathbb{R}^{n \times n}$ satisfying dual integral equation (\ref{G}).
    Moreover, the equality below is valid:
    \begin{equation} \label{F_G}
        G(t, s)
        = F(t, s),
        \quad (t, s) \in \Omega.
    \end{equation}
\end{prop}
\proof
    For any $t \in [t_0, \vartheta],$ the existence and uniqueness of a continuous solution $[t_0, t] \ni s \mapsto G(t, s) \in \mathbb{R}^{n \times n}$ of integral equation (\ref{G}) can be proved by the scheme from Proposition~\ref{Proposition_Existence_and_Boundedness_F}.

    Let us consider the function
    \begin{equation} \label{widetilde_F}
        \widetilde{F}(t, s)
        = \frac{\Id}{\Gamma(\alpha)}
        + \frac{(t - s)^{1 - \alpha}}{\Gamma(\alpha)} \int_{s}^{t} \frac{F(t, \tau) A(\tau)}{(t - \tau)^{1 - \alpha} (\tau - s)^{1 - \alpha}} \, \rd \tau,
        \ (t, s) \in \Omega.
    \end{equation}
    Then, due to Proposition~\ref{Proposition_Existence_and_Boundedness_F}, in order to prove equality (\ref{F_G}), it is sufficient to verify that, for any $s \in [t_0, \vartheta],$ the function $[s, \vartheta] \ni t \rightarrow \widetilde{F}(t, s) \in \mathbb{R}^{n \times n}$ is a continuous solution of integral equation (\ref{F}).
    For $s = \vartheta,$ this statement follows directly from relations (\ref{F}) and (\ref{widetilde_F}).

    Let $s \in [t_0, \vartheta)$ be fixed.
    We have
    \begin{equation*}
        \widetilde{F}(t, s)
        = \frac{\Id}{\Gamma(\alpha)} + (\J^\alpha_{s +} \varphi^{(t)})(t),
        \quad t \in [s, \vartheta],
    \end{equation*}
    where $\varphi^{(t)}(\tau) = F(t, \tau) A(\tau),$ $\tau \in [s, t].$
    According to the first inclusion in (\ref{A_b}) and Propositions~\ref{Proposition_Existence_and_Boundedness_F} and~\ref{Proposition_F_Holder}, we have $\varphi^{(t)}(\cdot) \in \Linf([s, t], \mathbb{R}^{n \times n})$ and $\|\varphi^{(t)}(\cdot)\|_{[s, t]} \leq M_A M_F,$ where $M_A$ is the number from (\ref{M_A}).
    Therefore, by Lemma~\ref{Lemma_J}, the function $[s, \vartheta] \ni t \mapsto \widetilde{F}(t, s) \in \mathbb{R}^{n \times n}$ is correctly defined, and the equality below is valid:
    \begin{equation*}
        \widetilde{F}(t, s)
        = \frac{\Id}{\Gamma(\alpha)}
        + (\I^\alpha_{s +} \psi^{(t)})(t),
        \quad t \in [s, \vartheta],
    \end{equation*}
    where $\psi^{(t)}(\tau) = \varphi^{(t)}(\tau) + (\R_{s +}^{\, \alpha} \varphi^{(t)})(\tau),$ $\tau \in [s, t].$
    Proposition~\ref{Proposition_R} yields
    \begin{multline*}
        \|\psi^{(t)}(\cdot)\|_{[s, t]}
        \leq \|\varphi^{(t)}(\cdot)\|_{[s, t]} + M_\R \|\varphi^{(t)}(\cdot)\|_{[s, t]} \\
        \leq (1 + M_\R) M_A M_F,
        \quad t \in [s, \vartheta].
    \end{multline*}
    Further, let $t_1, t_2 \in [s, \vartheta],$ $t_1 < t_2,$ be fixed.
    Since, due to Proposition \ref{Proposition_F_Holder},
    \begin{multline*}
        \|\varphi^{(t_1)}(\tau) - \varphi^{(t_2)}(\tau)\|
        = \| F(t_1, \tau) A(\tau) - F(t_2, \tau) A(\tau)\| \\
        \leq H_F M_A (t_2 - t_1)^\alpha
        \text{ for a.e. } \tau \in [s, t_1],
    \end{multline*}
    then, applying Proposition~\ref{Proposition_R} again, we derive
    \begin{multline*}
        \|\psi^{(t_1)}(\tau) - \psi^{(t_2)}(\tau)\| \\
        \leq \|\varphi^{(t_1)}(\tau) - \varphi^{(t_2)}(\tau)\|
        + \|(\R_{s +}^{\, \alpha} \varphi^{(t_1)})(\tau) - (\R_{s +}^{\, \alpha} \varphi^{(t_2)})(\tau)\| \\
        \leq (1 + M_\R) H_F M_A (t_2 - t_1)^\alpha
        \text{ for a.e. } \tau \in [s, t_1].
    \end{multline*}
    Thus, taking Proposition~\ref{Prop_I} into account, we conclude
    \begin{multline*}
        \|\widetilde{F}(t_1, s) - \widetilde{F}(t_2, s)\|
        = \|(\I^\alpha_{s +} \psi^{(t_1)})(t_1) - (\I^\alpha_{s +} \psi^{(t_2)})(t_2)\| \\
        \leq \|(\I^\alpha_{s +} \psi^{(t_1)})(t_1) - (\I^\alpha_{s +} \psi^{(t_2)})(t_1)\|
        + \|(\I^\alpha_{s +} \psi^{(t_2)})(t_1) - (\I^\alpha_{s +} \psi^{(t_2)})(t_2)\| \\
        \leq (1 + M_\R) H_\I H_F M_A (t_2 - t_1)^\alpha (t_1 - s)^\alpha
        + (1 + M_\R) H_\I M_A M_F (t_2 - t_1)^\alpha \\
        \leq ( (\vartheta - t_0)^\alpha H_F + M_F)(1 + M_\R) H_\I M_A (t_2 - t_1)^\alpha.
    \end{multline*}
    Hence, the function $[s, \vartheta] \ni t \mapsto \widetilde{F}(t, s) \in \mathbb{R}^{n \times n}$ is continuous, and it remains to verify that this function satisfies integral equation (\ref{F}).
    In order to do this, it is sufficient to show that, for $t \in [s, \vartheta]$,
    \begin{equation} \label{Prop_G_H_proof_equality}
        \int_{s}^{t} \frac{F(t, \tau) A(\tau)}{(t - \tau)^{1 - \alpha} (\tau - s)^{1 - \alpha}} \, \rd \tau
        = \int_{s}^{t} \frac{A(\tau) \widetilde{F}(\tau, s)}{(t - \tau)^{1 - \alpha} (\tau - s)^{1 - \alpha}} \, \rd \tau.
    \end{equation}
    For $t = s,$ equality (\ref{Prop_G_H_proof_equality}) holds automatically.
    Let $t \in (s, \vartheta]$ be fixed.
    In accordance with definition (\ref{F}) of the function $F,$ we obtain
    \begin{multline} \label{Prop_G_H_proof_equality_1}
        \int_{s}^{t} \frac{F(t, \tau) A(\tau)}{(t - \tau)^{1 - \alpha} (\tau - s)^{1 - \alpha}} \, \rd \tau
        = \frac{1}{\Gamma(\alpha)} \int_{s}^{t} \frac{A(\tau)}{(t - \tau)^{1 - \alpha} (\tau - s)^{1 - \alpha}} \, \rd \tau \\
        + \frac{1}{\Gamma(\alpha)} \int_{s}^{t} \int_{\tau}^{t} \frac{A(\xi) F(\xi, \tau) A(\tau)}
        {(\tau - s)^{1 - \alpha} (t - \xi)^{1 - \alpha} (\xi - \tau)^{1 - \alpha}} \, \rd \xi \, \rd \tau.
    \end{multline}
    Due to the first inclusion in (\ref{A_b}) and Proposition~\ref{Proposition_F_Holder}, by the Fubini--Tonelli theorem, we interchange the order of integration in the repeated integral:
    \begin{multline} \label{Prop_G_H_proof_equality_2}
        \frac{1}{\Gamma(\alpha)} \int_{s}^{t} \int_{\tau}^{t} \frac{A(\xi) F(\xi, \tau) A(\tau)}
        {(\tau - s)^{1 - \alpha} (t - \xi)^{1 - \alpha} (\xi - \tau)^{1 - \alpha}} \, \rd \xi \, \rd \tau \\
        = \frac{1}{\Gamma(\alpha)} \int_{s}^{t} \frac{A(\xi)}{(t - \xi)^{1 - \alpha}}
        \int_{s}^{\xi} \frac{F(\xi, \tau) A(\tau)}{(\tau - s)^{1 - \alpha} (\xi - \tau)^{1 - \alpha}} \, \rd \tau \, \rd \xi.
    \end{multline}
    From relations (\ref{Prop_G_H_proof_equality_1}) and (\ref{Prop_G_H_proof_equality_2}), taking definition (\ref{widetilde_F}) of the function $\widetilde{F}$ into account, we derive equality (\ref{Prop_G_H_proof_equality}).
    This completes the proof.
\proofend

Let us note that the use of the function $F$ instead of $Z$ may be more convenient in applications.
For example, it makes all the singularities in representation formulas (see {\it e.g.} relations (\ref{Cauchy_formula_PC}), (\ref{Cauchy_formula_GC}), and (\ref{Cauchy_formula_GC_compact}) below) readily seen.
Furthermore, numerical computing the values of the function $F$ seems simpler than of $Z,$ since $F$ does not have any singularities.

\section{Representation formulas}
\label{Section_Representation_formulas}

\setcounter{section}{4}
\setcounter{equation}{0}\setcounter{theorem}{0}

In this section, representation formulas for the solution $x(\cdot)$ of Cauchy problem (\ref{differential_equation}), (\ref{initial_condition}) are given.
First, following \cite[Theorem~5]{Bourdin_2018}, a particular case of initial condition (\ref{initial_condition}) when $t_\ast = t_0$ is studied.
After that, the general case is considered.

\begin{theo} \label{Theorem_Cauchy_formula_PC}
    Let $w_0 \in \mathbb{R}^n,$ and let $x(\cdot)$ be the solution of Cauchy problem (\ref{differential_equation}), (\ref{initial_condition}) for $t_\ast = t_0$ and $w_\ast(t_0) = w_0.$
    Then, the following representation formula holds:
    \begin{equation} \label{Cauchy_formula_PC}
        x(t)
        = \Big( \Id + \int_{t_0}^{t} \frac{F(t, \tau) A(\tau)}{(t - \tau)^{1 - \alpha}} \, \rd \tau \Big) w_0
        + \int_{t_0}^{t} \frac{F(t, \tau) b(\tau)}{(t - \tau)^{1 - \alpha}} \, \rd \tau,
    \end{equation}
    where $t \in [t_0, \vartheta],$ and the function $F$ is defined by (\ref{F}).
\end{theo}
\proof
    Let $w_0 \in \mathbb{R}^n$ be fixed.
    By the right-hand side of equality (\ref{Cauchy_formula_PC}), let us define the function
    \begin{equation} \label{Theorem_Cauchy_formula_PC_p_y}
        y(t)
        = w_0 + \int_{t_0}^{t} \frac{F(t, \tau) (A(\tau) w_0 + b(\tau))}{(t - \tau)^{1 - \alpha}} \, \rd \tau,
        \quad t \in [t_0, \vartheta].
    \end{equation}
    Based on inclusions (\ref{A_b}) and Proposition~\ref{Proposition_F_Holder}, one can obtain that the function $y(\cdot)$ is correctly defined, and, moreover, $y(\cdot) \in \C([t_0, \vartheta], \mathbb{R}^n).$
    Therefore, by Proposition~\ref{Proposition_x}, in order to prove equality (\ref{Cauchy_formula_PC}), it is sufficient to show that the function $y(\cdot)$ satisfies integral equation (\ref{integral_equation}), which in the considered particular case takes the following form:
    \begin{equation*}
        y(t)
        = w_0 + \frac{1}{\Gamma(\alpha)} \int_{t_0}^{t} \frac{A(\tau) y(\tau) + b(\tau)}{(t - \tau)^{1 - \alpha}} \, \rd \tau,
        \quad t \in [t_0, \vartheta].
    \end{equation*}
    Thus, it remains to verify that, for every $t \in [t_0, \vartheta],$
    \begin{equation} \label{Theorem_Cauchy_formula_PC_p_integral_equation}
        \int_{t_0}^{t} \frac{F(t, \tau) (A(\tau) w_0 + b(\tau))}{(t - \tau)^{1 - \alpha}} \, \rd \tau
        = \frac{1}{\Gamma(\alpha)} \int_{t_0}^{t} \frac{A(\tau) y(\tau) + b(\tau)}{(t - \tau)^{1 - \alpha}} \, \rd \tau.
    \end{equation}
    For $t = t_0,$ equality (\ref{Theorem_Cauchy_formula_PC_p_integral_equation}) holds automatically.
    Let $t \in (t_0, \vartheta]$ be fixed.
    According to definition (\ref{Theorem_Cauchy_formula_PC_p_y}) of the function $y(\cdot),$ we have
    \begin{multline} \label{Theorem_Cauchy_formula_PC_p_integral_equation_1}
        \int_{t_0}^{t} \frac{A(\tau) y(\tau) + b(\tau)}{(t - \tau)^{1 - \alpha}} \, \rd \tau
        = \int_{t_0}^{t} \frac{A(\tau) w_0 + b(\tau)}{(t - \tau)^{1 - \alpha}} \, \rd \tau \\
        + \int_{t_0}^{t} \int_{t_0}^{\tau} \frac{A(\tau) F(\tau, \xi) (A(\xi) w_0 + b(\xi))}
        {(t - \tau)^{1 - \alpha} (\tau - \xi)^{1 - \alpha}} \, \rd \xi \, \rd \tau.
    \end{multline}
    Due to inclusions (\ref{A_b}) and Proposition~\ref{Proposition_F_Holder}, by the Fubini--Tonelli theorem, one can interchange the order of integration in the repeated integral:
    \begin{multline} \label{Theorem_Cauchy_formula_PC_p_integral_equation_2}
        \int_{t_0}^{t} \int_{t_0}^{\tau} \frac{A(\tau) F(\tau, \xi) (A(\xi) w_0 + b(\xi))}
        {(t - \tau)^{1 - \alpha} (\tau - \xi)^{1 - \alpha}} \, \rd \xi \, \rd \tau \\
        = \int_{t_0}^{t} \int_{\xi}^{t} \frac{A(\tau) F(\tau, \xi)}{(t - \tau)^{1 - \alpha} (\tau - \xi)^{1 - \alpha}} \, \rd \tau
        \, (A(\xi) w_0 + b(\xi)) \, \rd \xi.
    \end{multline}
    From relations (\ref{Theorem_Cauchy_formula_PC_p_integral_equation_1}) and (\ref{Theorem_Cauchy_formula_PC_p_integral_equation_2}), taking definition (\ref{F}) of the function $F$ into account, we derive equality (\ref{Theorem_Cauchy_formula_PC_p_integral_equation}).
    This completes the proof.
\proofend

An extension of representation formula (\ref{Cauchy_formula_PC}) to the general case of initial condition (\ref{initial_condition}) is based on the following technical lemma.
\begin{lem} \label{Lem_Continuation}
    Let $t_\ast \in (t_0, \vartheta),$ $\varphi(\cdot) \in \Linf([t_0, t_\ast], \mathbb{R}^n),$ and
    \begin{equation} \label{Lem_Continuation_psi}
        \psi(t)
        = \frac{\sin (\alpha \pi)}{\pi} \int_{t_0}^{t_\ast} \frac{(t_\ast - \tau)^\alpha \varphi(\tau)}{t - \tau} \, \rd \tau,
        \quad t \in [t_\ast, \vartheta].
    \end{equation}
    Then, the inclusion $\psi(\cdot) \in \C([t_\ast, \vartheta], \mathbb{R}^n)$ is valid, and the equalities
    \begin{equation} \label{Lem_Continuation_equalities}
        \begin{array}{c}
            \displaystyle
            \int_{t_0}^{t_\ast} \frac{\varphi(\tau)}{(t - \tau)^{1 - \alpha}} \, \rd \tau
            = \int_{t_\ast}^{t} \frac{\psi(\tau)}{(t - \tau)^{1 - \alpha} (\tau - t_\ast)^\alpha} \, \rd \tau, \\[1em]
            \displaystyle
            \psi(t)
            = \frac{(t - t_\ast)^\alpha \alpha}{\Gamma(1 - \alpha)} \int_{t_0}^{t_\ast} \frac{x(\tau)}{(t - \tau)^{1 + \alpha}} \, \rd \tau
        \end{array}
    \end{equation}
    hold for $t \in (t_\ast, \vartheta],$ where $x(t) = (\I^\alpha_{t_0 +} \varphi)(t),$ $t \in [t_0, t_\ast].$
\end{lem}
\proof
    The validity of the inclusion $\psi(\cdot) \in \C([t_\ast, \vartheta], \mathbb{R}^n)$ follows directly from relation (\ref{Lem_Continuation_psi}).
    Equalities (\ref{Lem_Continuation_equalities}) can be proved by the scheme from \cite[Theorem~13.10]{Samko_Kilbas_Marichev_1993}.
\proofend

\begin{theo} \label{Theorem_Cauchy_formula_GC}
    Let $t_\ast \in [t_0, \vartheta),$ $w_\ast(\cdot) \in \AC([t_0, t_\ast], \mathbb{R}^n),$ and let $x(\cdot)$ be the solution of Cauchy problem (\ref{differential_equation}), (\ref{initial_condition}).
    Then, the following representation formula holds:
    \begin{equation} \label{Cauchy_formula_GC}
        x(t) = \Big( \Id + \int_{t_\ast}^{t} \frac{F(t, \tau) A(\tau)}{(t - \tau)^{1 - \alpha}} \, \rd \tau \Big) w_\ast(t_\ast)
        + \int_{t_\ast}^{t} \frac{F(t, \tau) b_\ast(\tau)}{(t - \tau)^{1 - \alpha}} \, \rd \tau,
    \end{equation}
    where $t \in [t_\ast, \vartheta],$ the function $F$ is defined by (\ref{F}), and, for $t \in (t_\ast, \vartheta],$
    \begin{equation} \label{b_ast}
        b_\ast(t)
        = \frac{\alpha}{\Gamma(1 - \alpha)} \int_{t_0}^{t_\ast} \frac{w_\ast(\tau) - w_\ast(t_0)}{(t - \tau)^{1 + \alpha}} \, \rd \tau
        - \frac{w_\ast(t_\ast) - w_\ast(t_0)}{\Gamma(1 - \alpha) (t - t_\ast)^\alpha}
        + b(t).
    \end{equation}
\end{theo}
\proof
    Let $t_\ast \in [t_0, \vartheta)$ and $w_\ast(\cdot) \in \AC([t_0, t_\ast], \mathbb{R}^n)$ be fixed.
    In the case when $t_\ast = t_0,$ we have $b_\ast(t) = b(t),$ $t \in (t_0, \vartheta],$ and formula (\ref{Cauchy_formula_GC}) follows from Theorem~\ref{Theorem_Cauchy_formula_PC}.

    Let $t_\ast > t_0.$
    By Proposition~\ref{Proposition_D}, let us choose $\varphi_\ast(\cdot) \in \Linf([t_0, t_\ast], \mathbb{R}^n)$ such that $\varphi_\ast(t) = (^C \D_{t_0 +}^{\, \alpha} w_\ast)(t)$ for almost every $t \in [t_0, t_\ast].$
    By this function $\varphi_\ast(\cdot),$ let us define $\psi_\ast(\cdot) \in \C([t_\ast, \vartheta], \mathbb{R}^n)$ according to Lemma~\ref{Lem_Continuation}.
    Hence,
    \begin{equation} \label{Theorem_Cauchy_formula_GC_p_psi_ast}
        \psi_\ast(t_\ast)
        = \frac{\sin (\alpha \pi)}{\pi} \int_{t_0}^{t_\ast} \frac{(^C \D_{t_0 +}^{\, \alpha} w_\ast)(\tau)}{(t_\ast - \tau)^{1 - \alpha}} \, \rd \tau
        = \frac{w_\ast(t_\ast) - w_\ast(t_0)}{\Gamma(1 - \alpha)},
    \end{equation}
    and, for any $t \in (t_\ast, \vartheta],$
    \begin{equation} \label{Theorem_Cauchy_formula_GC_p_equalities}
        \begin{array}{c}
            \displaystyle
            \int_{t_0}^{t_\ast} \frac{(^C \D_{t_0 +}^{\, \alpha} w_\ast)(\tau)}{(t - \tau)^{1 - \alpha}} \, \rd \tau
            = \int_{t_\ast}^{t} \frac{\psi_\ast(\tau)}{(t - \tau)^{1 - \alpha} (\tau - t_\ast)^\alpha} \, \rd \tau, \\[1em]
            \displaystyle
            \psi_\ast(t)
            = \frac{(t - t_\ast)^\alpha \alpha}{\Gamma(1 - \alpha)}
            \int_{t_0}^{t_\ast} \frac{w_\ast(\tau) - w_\ast(t_0)}{(t - \tau)^{1 + \alpha}} \, \rd \tau.
        \end{array}
    \end{equation}
    Thus, for the function $b_\ast(\cdot)$ from (\ref{b_ast}), we conclude
    \begin{equation} \label{Theorem_Cauchy_formula_GC_p_b_ast}
        b_\ast(t)
        = \frac{\psi_\ast(t) - \psi_\ast(t_\ast)}{(t - t_\ast)^\alpha}
        + b(t),
        \quad t \in (t_\ast, \vartheta].
    \end{equation}
    By initial condition (\ref{initial_condition}) and the right-hand side of equality (\ref{Cauchy_formula_GC}), let us define the function
    \begin{equation} \label{Theorem_Cauchy_formula_GC_p_y}
        y(t) =
        \left\{
            \begin{array}{l}
                w_\ast(t), \quad t \in [t_0, t_\ast], \\[1em]
                \displaystyle
                \Big( \Id + \int_{t_\ast}^{t} \frac{F(t, \tau) A(\tau)}{(t - \tau)^{1 - \alpha}} \, \rd \tau \Big) w_\ast(t_\ast) \\[1em]
                \displaystyle
                \hphantom{aaaaaaaaaaaaaaaa} + \int_{t_\ast}^{t} \frac{F(t, \tau) b_\ast(\tau)}{(t - \tau)^{1 - \alpha}} \, \rd \tau,
                \quad t \in (t_\ast, \vartheta].
            \end{array}
        \right.
    \end{equation}
    Based on the inclusions $\psi_\ast(\cdot) \in \C([t_\ast, \vartheta], \mathbb{R}^n)$ and (\ref{A_b}) and Proposition~\ref{Proposition_F_Holder}, one can obtain that the function $y(\cdot)$ is correctly defined, and, moreover, $y(\cdot) \in \C([t_0, \vartheta], \mathbb{R}^n).$
    Therefore, owing to Proposition~\ref{Proposition_x}, in order to prove equality (\ref{Cauchy_formula_GC}), it is sufficient to show that $y(\cdot)$ satisfies the integral equation
    \begin{multline} \label{Theorem_Cauchy_formula_GC_p_integral_equation_basic}
        y(t) = w_\ast(t_0)
        + \frac{1}{\Gamma(\alpha)} \int_{t_0}^{t_\ast} \frac{(^C \D_{t_0 +}^{\, \alpha} w_\ast)(\tau)}{(t - \tau)^{1 - \alpha}} \, \rd \tau \\
        + \frac{1}{\Gamma(\alpha)} \int_{t_\ast}^{t} \frac{A(\tau) y(\tau) + b(\tau)}{(t - \tau)^{1 - \alpha}} \, \rd \tau,
        \quad t \in [t_\ast, \vartheta].
    \end{multline}
    For $t = t_\ast,$ equality (\ref{Theorem_Cauchy_formula_GC_p_integral_equation_basic}) follows from the relations
    \begin{equation*}
        y(t_\ast)
        = w_\ast(t_\ast)
        = w_\ast(t_0) + \frac{1}{\Gamma(\alpha)} \int_{t_0}^{t_\ast} \frac{(^C \D_{t_0 +}^{\, \alpha} w_\ast)(\tau)}{(t_\ast - \tau)^{1 - \alpha}} \, \rd \tau,
    \end{equation*}
    which are valid due to definition (\ref{Theorem_Cauchy_formula_GC_p_y}) of the function $y(\cdot)$ and Proposition~\ref{Proposition_D}.
    Let $t \in (t_\ast, \vartheta]$ be fixed.
    By relations (\ref{Theorem_Cauchy_formula_GC_p_psi_ast}) and the first equality in (\ref{Theorem_Cauchy_formula_GC_p_equalities}), for the first two terms in (\ref{Theorem_Cauchy_formula_GC_p_integral_equation_basic}), we derive
    \begin{multline*}
        w_\ast(t_0)
        + \frac{1}{\Gamma(\alpha)} \int_{t_0}^{t_\ast} \frac{(^C \D_{t_0 +}^{\, \alpha} w_\ast)(\tau)}{(t - \tau)^{1 - \alpha}} \, \rd \tau \\
        = w_\ast(t_\ast)
        + \frac{1}{\Gamma(\alpha)} \int_{t_\ast}^{t} \frac{\psi_\ast(\tau) - \psi_\ast(t_\ast)}{(t - \tau)^{1 - \alpha} (\tau - t_\ast)^\alpha} \, \rd \tau.
    \end{multline*}
    Thus, according to equality (\ref{Theorem_Cauchy_formula_GC_p_b_ast}), in order to complete the proof, it remains to verify that
    \begin{equation*}
        \int_{t_\ast}^{t} \frac{F(t, \tau) (A(\tau) w_\ast(t_\ast) + b_\ast(\tau))}{(t - \tau)^{1 - \alpha}} \, \rd \tau
        = \frac{1}{\Gamma(\alpha)} \int_{t_\ast}^{t} \frac{A(\tau) y(\tau) + b_\ast(\tau)}{(t - \tau)^{1 - \alpha}} \, \rd \tau.
    \end{equation*}
    This equality is proved by the same steps as equality (\ref{Theorem_Cauchy_formula_PC_p_integral_equation}).
    The only difference is that the function $b_\ast(\cdot)$ is not essentially bounded in general, and, therefore, when applying the Fubini--Tonelli theorem, equality (\ref{Theorem_Cauchy_formula_GC_p_b_ast}) and the inclusion $\psi_\ast(\cdot) \in \C([t_\ast, \vartheta], \mathbb{R}^n)$ should be used.
\proofend

In the corollary below, representation formula (\ref{Cauchy_formula_GC}) is rewritten in another form.
\begin{cor}
    Let $t_\ast \in [t_0, \vartheta),$ $w_\ast(\cdot) \in \AC([t_0, t_\ast], \mathbb{R}^n),$ and let $x(\cdot)$ be the solution of Cauchy problem (\ref{differential_equation}), (\ref{initial_condition}).
    Then, the following representation formula holds:
    \begin{multline} \label{Cauchy_formula_GC_compact}
            x(t) =
            \Big( \Id + \int_{t_\ast}^{t} \frac{F(t, \tau) A(\tau)}{(t - \tau)^{1 - \alpha}} \, \rd \tau \Big) w_\ast(t_0) \\
            + \frac{\alpha}{\Gamma(1 - \alpha)} \int_{t_\ast}^{t} \frac{F(t, \tau)}{(t - \tau)^{1 - \alpha}}
            \int_{t_0}^{t_\ast} \frac{w_\ast(\xi) - w_\ast(t_0)}{(\tau - \xi)^{1 + \alpha}} \, \rd \xi \, \rd \tau \\
            + \int_{t_\ast}^{t} \frac{F(t, \tau) b(\tau)}{(t - \tau)^{1 - \alpha}} \, \rd \tau,
            \quad t \in (t_\ast, \vartheta].
    \end{multline}
\end{cor}
\proof
    Let $t_\ast \in [t_0, \vartheta)$ and $w_\ast(\cdot) \in \AC([t_0, t_\ast], \mathbb{R}^n)$ be fixed.
    If $t_\ast = t_0,$ then formula (\ref{Cauchy_formula_GC_compact}) follows from Theorem~\ref{Theorem_Cauchy_formula_PC}, since the second term in the right-hand side of this formula vanishes.

    In order to prove equality (\ref{Cauchy_formula_GC_compact}) in the case when $t_\ast > t_0,$ in accordance with Theorem~\ref{Theorem_Cauchy_formula_GC}, it is sufficient to show that, for every $t \in (t_\ast, \vartheta],$
    \begin{equation} \label{Corollary_Cauchy_formula_GC_compact_p_equality}
        \Id + \int_{t_\ast}^{t} \frac{F(t, \tau) A(\tau)}{(t - \tau)^{1 - \alpha}} \, \rd \tau
        = \frac{1}{\Gamma(1 - \alpha)} \int_{t_\ast}^{t} \frac{F(t, \tau)}{(t - \tau)^{1 - \alpha} (\tau - t_\ast)^\alpha} \, \rd \tau.
    \end{equation}
    Let $t \in (t_\ast, \vartheta]$ be fixed.
    By Proposition~\ref{Proposition_F_G}, we derive
    \begin{multline} \label{Corollary_Cauchy_formula_GC_compact_p_equality_1}
        \frac{1}{\Gamma(1 - \alpha)} \int_{t_\ast}^{t} \frac{F(t, \tau)}{(t - \tau)^{1 - \alpha} (\tau - t_\ast)^\alpha} \, \rd \tau
        = \Id \\
        + \frac{1}{\Gamma(1 - \alpha) \Gamma(\alpha)} \int_{t_\ast}^{t} \int_{\tau}^{t} \frac{F(t, \xi) A(\xi)}
        {(\tau - t_\ast)^\alpha (t - \xi)^{1 - \alpha} (\xi - \tau)^{1 - \alpha}} \, \rd \xi \, \rd \tau.
    \end{multline}
    Due to the first inclusion in (\ref{A_b}) and Proposition~\ref{Proposition_F_Holder}, by the Fubini--Tonelli theorem, we interchange the order of integration in the repeated integral:
    \begin{multline} \label{Corollary_Cauchy_formula_GC_compact_p_equality_2}
        \frac{1}{\Gamma(1 - \alpha) \Gamma(\alpha)} \int_{t_\ast}^{t} \int_{\tau}^{t} \frac{F(t, \xi) A(\xi)}
        {(\tau - t_\ast)^\alpha (t - \xi)^{1 - \alpha} (\xi - \tau)^{1 - \alpha}} \, \rd \xi \, \rd \tau \\
        = \int_{t_\ast}^{t} \frac{F(t, \xi) A(\xi)}{(t - \xi)^{1 - \alpha}} \, \rd \xi.
    \end{multline}
    Equality (\ref{Corollary_Cauchy_formula_GC_compact_p_equality}) follows from relations (\ref{Corollary_Cauchy_formula_GC_compact_p_equality_1}) and (\ref{Corollary_Cauchy_formula_GC_compact_p_equality_2}).
\proofend

It seems that representation formula (\ref{Cauchy_formula_GC_compact}) is simpler than formula (\ref{Cauchy_formula_GC}).
However, it should be noted that, in general, the second term in the right-hand side of formula (\ref{Cauchy_formula_GC_compact}) does not tend to zero when $t$ goes to $t_\ast.$
This makes it necessary to use this formula more carefully, especially when developing the corresponding numerical methods.

\section*{Acknowledgements}

This work was supported by RSF, project no. 19-11-00105.

 \bigskip \smallskip

 \it

 \noindent
Krasovskii Institute of Mathematics and Mechanics \\
Ural Branch of Russian Academy of Sciences \\
S. Kovalevskaya Str., Block 16 \\
Ekaterinburg -- 620990, RUSSIA  \\[4pt]
Ural Federal University \\
Mira Str., Block 32 \\
Ekaterinburg -- 620002, RUSSIA  \\[4pt]
e-mail: m.i.gomoyunov@gmail.com
\hfill Received: August 22, 2019 \\[12pt]

\end{document}